\documentclass[a4,11pt]{article}
\usepackage[sumlimits]{amsmath}
\usepackage{amssymb,amscd}
\usepackage[francais]{babel}

\usepackage{graphicx}  

\def\nal{\| \hspace{-.10em} |}

\def\1{{\mathbf 1 }}

\catcode`\é=\active\def é{\'e} \catcode`\è=\active\def è{\`e}
\catcode`\ç=\active\def ç{\c{c}} \catcode`\à=\active\def à{\`a}
\catcode`\ê=\active\def ê{\^e} \catcode`\ù=\active\def ù{\`u}
\catcode`\ô=\active\def ô{\^o} \catcode`\û=\active\def û{\^u}
\catcode`\î=\active\def î{\^{\i}} \catcode`\â=\active\def â{\^a}
\catcode`\ö=\active\def ö{\"o}
\catcode`\ë=\active\def ë{\"e}

\def\nnn{{\cal N}}

\def\ph{\varphi}
\def\eps{{\varepsilon}}

\def\Om{{\Omega}}

\def\VaO{V_{\alpha}(\Omega)}

\def \Nal{\mathcal N}

\def \Ual{\mathcal U}
\def \Val{\mathcal V}

\newtheorem{theo}{Th\'eor\`eme}

\newtheorem{prop}[theo]{Proposition}

\newtheorem{exercice}{\bf Exercice}

\def\N{{\mathbb N}}
\def\R{{\mathbb R}}
\def\Z{{\mathbb Z}}

\def\Om{{\Omega}}

\def\ph{\varphi}
\def\eps{{\varepsilon}}

\def\cov{{\rm Cov}}

\def\1{{\mathbf 1}}

\title{Décroissance exponentielle des corrélations pour un système dynamique réel induit d'un système en dimension 2\\
}
\author{JAGER Lisette, MAES Jules, NINET Alain}

\begin{document}

\maketitle

%
%

\section*{Résumé}

\noindent On étudie le processus $\{X_n, n\in \N \}$ à variables réelles bornées défini par $X_{n+2} = \varphi(X_n,X_{n+1})$. On cherche à obtenir la décroissance des corrélations de ce modèle en prenant uniquement des hypothèses analytiques sur la transformation.

%
%

\section{Introduction}

\noindent Depuis les années 80, l'étude des séries temporelles non-linéaires par les statisticiens a permis  de modéliser de nombreux phénomènes en Physique, \'Economie et en Finance. Mais c'est dans les années 90 que  la théorie du Chaos est devenue un axe de recherche incontournable pour l'étude de ces processus. Pour une revue complète de cette littérature, sur la théorie du Chaos on pourra consulter Collet-Eckmann \cite{CE} et concernant l'étude des séries temporelles non-linéaires Chan-Tong \cite{TON1}. Dans cette perspective, un modèle général pourrait s'écrire:
$$  X_{t+1}= \ph(X_{t},\dots, X_{t-d+1}) + \varepsilon_t,$$
avec $\varphi $ non-linéaire et $ \varepsilon_t$ un bruit. Nous proposons une première étude du "squelette" de ce modèle comme l'appelle Tong \cite{TON2} en commençant avec $d=2$, c'est à dire par une étude d'un système dynamique induit de ce modèle. En effet, on considère le modèle à variables bornées, $X_{n+2}= \ph(X_n,X_{n+1})$, avec $\varphi :\Ual^2 \rightarrow \Ual$ pour $\Ual=[-L,L]$ et $L \in \R_+^*$, $\varphi$ étant définie par morceaux sur $\Ual^2$. Ce modèle induit un système dynamique  $ (\Om,\tau ,\mu,T)$ où $\mu$ est une mesure invariante par la transformation $T :\Om \rightarrow \Om$, et $\Om$ un compact de $\R^2$.
Nous obtenons  sous certaines hypothèses sur $\varphi$, telles que $T$ vérifie les hypothèses de Saussol \cite{SAU}, la décroissance exponentielle des corrélations si $T$ est
 mélange. Plus précisément, pour des applications $f$ et $h$ bien choisies, il existe
 une constante $C=C(f,h)>0$, $ 0<\rho<1$ tels que :
 $$ \left|  \int_{\Om} f \circ T^n \,  h \ d\mu - \int_{\Om} f d\mu ~ \int_{\Om}h d\mu  \right| \leqslant C \, \rho^n.  $$
Ce résultat nous permet d'en déduire une inégalité de covariance du type : 
$$ \left| \, \cov ( \, f(X_n),h(X_0) \,) \,  \right| \leqslant C \, \rho^n .$$ 
D'autres voies seraient possibles pour obtenir le même résultat sous différentes hypothèses sur le système induit. On peut en effet citer la méthode des "tours de Young" \cite{YOU}. Pour une vue d'ensemble de ces techniques, on pourra consulter l'article de Alves-Freitas-Luzzato-Vaienti \cite{AFLV}.\\
A la fin de cet article, nous  illustrons nos résultats par deux exemples, l'un linéaire par morceaux et l'autre non linéaire.\\

%
%

\section{Hypothèses et résultats}

\noindent Soit $\varphi : [-L,L]^2 \rightarrow [-L,L]$ pour un certain $L \in \R_+^*$, définie par morceaux sur $[-L,L]^2$. Pour étudier le processus $\{X_n, n\in \N \}$ défini par $X_{n+2}= \ph(X_n,X_{n+1})$, il
existe différentes façons de choisir le système dynamique induit $Z_{n+1} = T(Z_n)$ avec $Z_n \in \R^2$.
Nous avons tenté deux approches, d'une part la méthode canonique, en prenant
 $T(x,y) = (y, \varphi(x,y))$, et d'autre part une double itération, ce qui revient à prendre $T(x,y) = (\varphi(x,y),\varphi(y,\varphi(x,y)))$.
Il s'avère que la première méthode, moyennant une conjugaison ($T(x,y) = (\frac{y}{\gamma},\gamma \varphi(x,\frac{y}{\gamma}))$ avec $Z_n = (X_n,\gamma X_{n+1})$) est la plus fructueuse, les hypothèses à imposer pour faire fonctionner la seconde approche étant plus lourdes et les résultats obtenus moins bons. Il a alors été possible de travailler dans des espaces similaires aux espaces $V_{\alpha}$ de Saussol et d'utiliser ses résultats.\\

\noindent Plus précisément, on suppose que les hypothèses suivantes sont vérifiées :
\begin{enumerate}
\item 
Il existe $d \in \N^*$ tel que
$$
[-L,L]^2= \bigcup_{k=1}^d O_k \ \cup \Nal,
$$
où les $O_k$ sont des ouverts non vides, $\Nal$ est de mesure de Lebesgue nulle et la réunion est disjointe. Les bords des $O_k$ sont inclus dans des sous-variétés compactes $C^1$ de dimension $1$ de $\R^2$. 
\item
Il existe $\eps_1>0$ tel que, pour tout $k\in \{1,\dots d\}$, il existe une application $\ph_k$  définie sur $B_{\eps_1}(\overline{O_k}) = \{ (x,y) \in \R^2, ~ d((x,y),\overline{O_k}) \leq \eps_1\}$ à valeurs dans $\R$, telle que $\ph_k|_{O_k} = \ph|_{O_k}$.
\item
L'application $\ph_k$ est  bornée, de classe $C^{1,\alpha}$ sur $B_{\eps_1}(\overline{O_k})$ pour un $\alpha\in ]0,1]$ \footnote{Si $\ph_k$ est  de classe $C^2$ sur $B_{\eps_1}(\overline{O_k})$, elle est nécessairement de classe $C^{1,\alpha}$ sur $B_{\eps_1}(\overline{O_k})$ avec $\alpha = 1$}, ce qui signifie que $\ph_k$ est de classe $C^1$ et qu'il existe $C_k>0$ tel que, pour tous $ (u,v), (u',v')$ de $B_{\eps_1}(\overline{O_k})$, 
$$
\begin{array}{lll}
\displaystyle 
\left|\frac{\partial \ph_k}{\partial u}(u,v) -
\frac{\partial \ph_k}{\partial u}(u',v')\right|
\leq C_k ||(u,v)-(u',v')||^{\alpha}\\
\displaystyle 
\left|\frac{\partial \ph_k}{\partial v}(u,v) -
\frac{\partial \ph_k}{\partial v}(u',v')\right|
\leq C_k ||(u,v)-(u',v')||^{\alpha}.\\
\end{array}
$$
On suppose également qu'il existe $A>1$ et $M\in ]0,A-1[$ tels que  :
$$ 
 \forall (u,v) \in B_{\eps_1}(\overline{O_k}),\qquad
 \left|\frac{\partial \ph_k}{\partial u}(u,v)
\right|\geq A,\quad
\left|\frac{\partial \ph_k}{\partial v}(u,v)
\right| \leq M.
$$
\item 
Les $O_k$ vérifient la condition géométrique qui suit \footnote{Dans les cas favorables, l'hypothèse géométrique peut être remplacée par la suivante, plus forte mais plus simple : pour tous points $(u,v)$ et $(u',v)$ de $B_{\eps_1}(\overline{O_k})$, le segment $[(u,v),(u',v)]$ est inclus dans   $B_{\eps_1}(\overline{O_k})$} : pour tous  $(u,v)$ et $(u',v)$ de $ B_{\eps_1}(\overline{O_k})$, il existe un chemin $\Gamma=(\Gamma_1,\Gamma_2) : [0,1]\rightarrow  B_{\eps_1}(\overline{O_k})$, de classe $C^1$, reliant $(u,v)$ à $(u',v)$, de gradient non nul et vérifiant
\begin{equation*}
\forall t\in ]0,1[, \left|\Gamma_1'(t)\right|
 >\frac{M}{A}  \left|\Gamma_2'(t)\right| .
\end{equation*}
\item
Le nombre maximal d'arcs $C^1$ de $\Nal$ se croisant en un point est $Y \in \N^*$.
De plus, on pose 
$$
s = \left( \frac{  2A+M^2 -M\sqrt{M^2+4A} }{2} \right)^{-1/2} < 1
$$
et l'on suppose que 
$$
\eta:=  s^{\alpha} + \frac{8 s}{\pi(1-s)}Y <1.
$$
\end{enumerate}

\noindent On notera $\displaystyle \gamma = \frac{1}{\sqrt{A}}<1$ et, pour tout $k \in \{1, ..., d\}$, $U_k$ (resp. $W_k$, $\Nal'$) sera l'image de $O_k$ (resp. $B_{\eps_1}(\overline{O_k})$, $\Nal$) par l'affinité qui à $(u,v) \in \R^2$ associe $(u,\gamma v)$. L'ensemble $\Om = [-L,L] \times [-\gamma L, \gamma L]$, sur lequel nous travaillerons, sera l'image de $[-L,L]^2$ par cette même affinité.\\
\\
Pour tout borélien de mesure non nulle $S$ de $\R^2$, pour toute $f \in L^1_m(\R^2,\R)$, on notera :
$$ Osc(f,S) = \underset{S}{Esup} f - \underset{S}{Einf} f $$
où $\underset{S}{Esup}$ et $\underset{S}{Einf}$ représentent le sup et l'inf essentiels sur $S$ pour la mesure de Lebesgue $m$.\\
On définit alors : 
$$
|f|_{\alpha}= \sup_{0<\eps<\eps_1}\eps^{-\alpha} \int_{\R^2}  {\rm Osc}(f,B_{\eps}(x,y))\ dxdy \qquad {\rm et} \qquad \| f \|_{\alpha} = \| f \|_{L^1_m} + |f|_{\alpha}
$$
puis l'ensemble $V_{\alpha} = \{ f \in L^1_m(\R^2,\R), ~ \| f \|_{\alpha} < +\infty \}$.\\

\noindent Nous introduisons des notions similaires sur $\Om$ : pour tout $0 < \eps_0 < \gamma \eps_1$, pour tout $g \in L^{\infty}_m(\Om,\R)$, on définit :
$$
N(g,\alpha,L) = \sup_{0<\eps<\eps_0}\eps^{-\alpha} \int_{\Om} {\rm Osc}(g,B_{\eps}(x,y) \cap \Om)\ dxdy.
$$
On pose alors :
$$
||g||_{\alpha,L}=  N(g,\alpha,L) + 16 (1+\gamma) \eps_0^{1-\alpha} L ||g||_{\infty} + ||g||_{L^1_m}.
$$
On dit que $g$ appartient à $\VaO$ si cette quantité est finie. Cet ensemble ne dépend pas du choix de $\eps_0$, mais $N$ et $\| . \|_{\alpha,L}$ en dépendent.\\
Il existe des liens entre ces deux ensembles. En effet, on peut démontrer le résultat suivant grâce à la proposition 3.4 de \cite{SAU} :

\begin{prop}
\begin{enumerate}
\item
Si $g \in \VaO$ et si on prolonge $g$ en une fonction notée $f$ en posant $f(x,y)=0$ si $(x,y)\notin \Om$, alors $f \in V_{\alpha}$ et 
$$
\| f \|_{\alpha} \leq \| g \|_{\alpha,L}.
$$
\item
Soit $f\in V_{\alpha}$. Posons $g= f \1_{\Om}$. Alors $g\in \VaO$ et l'on a 
$$
 \| g \|_{\alpha,L}\leq \left( 1+16 (1+\gamma) L
 \frac{\max(1,\eps_0^{\alpha})}{\pi \eps_0^{1+\alpha}}\right)
\| f \|_{\alpha}.
$$
\end{enumerate}
\end{prop}

\noindent On obtient ainsi, sous les hypothèses 1 à 5 précédentes, un premier résultat :

\begin{theo}\label{resconjug}
Soit $T$ la transformation définie sur $\Om$ par : $\forall  (x,y) \in U_k$ :
$$
T(x,y) = T_k(x,y)= \left( \frac{y}{\gamma},\  \gamma \ph_k(x,   \frac{y}{\gamma})\right).
$$
La définition des $T_k$ s'étend naturellement à $W_k$ par la même formule.
Alors :
\begin{enumerate}
\item 
L'opérateur de Perron-Frobénius $P : L^1_m(\Om) \rightarrow L^1_m(\Om)$ associé à $T$ a un nombre fini de valeurs  propres $\lambda_1,\dots,\lambda_r$ de module $1$.
\item 
Pour tout $i\in\{1,\dots,r\}$, l'espace propre $E_i=\{ f\in L^1_m(\Om) \ : \ Pf=\lambda_i f\}$  associé à la valeur propre $\lambda_i$ est inclus dans $V_{\alpha}(\Om)$ et de dimension finie.
\item 
L'opérateur $P$ se décompose en 
$$
P= \sum_{i=1}^r \lambda_i P_i + Q,
$$
où les $P_i$ sont des projections sur les espaces $E_i$, $\nal P_i \nal_1\leq 1$ et $Q$ est un opérateur linéaire sur $ L^1_m(\Om)$, vérifiant $Q(\VaO) \subset \VaO$, $\sup\limits_{n \in \N^*}\nal Q^n \nal_1<\infty$ et  $\nal Q^n \nal_{\alpha,L} = O(q^n)$ quand $n \rightarrow +\infty$ pour un $q\in ]0,1[$. De plus $P_iP_j=0$ si $i\neq j$, $P_iQ=QP_i=0$ pour tout $i$.
\item 
L'opérateur $P$ admet $1$ comme valeur propre. Supposons que $\lambda_1=1$, soit $h_*=P_1 \1_{\Om}$ et soit $d\mu=h_* ~ dm$. Alors $\mu$ est la plus grande mesure invariante absolument continue (ACIM) de $ T$, au sens suivant :
si $\nu <<m$ et $\nu$ est $T$ invariante, alors $\nu<<\mu$.
\item 
Le support de $\mu$ peut être décomposé en un nombre fini d'ensembles  mesurables, deux à deux disjoints, sur lesquels une puissance de $T$ est mélange (mixing). Plus précisément, pour tout $j \in \{1,2,\dots, \dim(E_1)\}$, il existe $L_j \in \N^*$ et $L_j$ ensembles $W_{j,l}$ $(0 \leq l \leq L_{j}-1)$ deux à deux disjoints, vérifiant $T(W_{j,l})=W_{j,l+1 \mod(L_j)}$ et $T^{L_j}$ est mélange sur chaque $W_{j,l}$.
On note $\mu_{j,l}$ la ``restriction normalisée`` de $\mu$ à $W_{j,l}$ définie par
$$
\mu_{j,l}(B)= \frac{\mu(B\cap W_{j,l})}{\mu(W_{j,l})}, \ d\mu_{j,l} = \frac{h^* \1_{W_{j,l}} }{\mu(W_{j,l})} dm.
$$
Dire que $T^{L_j}$ est mélange sur chaque $W_{j,l}$ signifie que, pour tout $f \in L^1_{\mu_{j,l}}(W_{j,l})$ et pour tout $h \in L^{\infty}_{\mu_{j,l}}(W_{j,l})$,
$$ \lim\limits_{n \rightarrow + \infty} <T^{nL_j} f,h>_{\mu_{j,l}} = <f,1>_{\mu_{j,l}} <1,h>_{\mu_{j,l}} $$
avec les notations indifféremment employées : $<f,g>_{\mu'} = \mu'(fg) = \int f g ~ d\mu'$. 
\item De plus, il existe $C>0$ et $0<\rho < 1$ tels que, pour tout $h$ dans $\VaO$ et $f\in L^1_{\mu}(\Om)$, on a 
$$
\left| \int_{\Om} f \circ T^{ n \times ppcm(L_i)}  h \ d\mu -\sum_{j=1}^{\dim(E_1)} \sum_{l=0}^{L_j-1} \mu(W_{j,l}) <f,1>_{\mu_{j,l}} <1,h>_{\mu_{j,l}} \right| \leq C  ||h||_{\alpha,\Om} ||f||_{L^1_{\mu}(\Om)}~\rho^{n }.
$$
\item Si de plus $T$ est mélange\footnote{Ce qui équivaut à : si 1 est la seule valeur propre de $P$ de module 1 et si de plus elle est simple}, alors le résultat précédent peut s'écrire :
 il existe $C>0$ et $0<\rho < 1$ tels que, pour tout $h$ dans $\VaO$ et $f\in L^1_{\mu}(\Om)$, on a  :
$$ \left|  \int_{\Om} f \circ T^n \,  h \ d\mu - \int_{\Om} f d\mu ~ \int_{\Om} h d\mu  \right| 
 \leq C  ||h||_{\alpha,\Om}~||f||_{L^1_{\mu}(\Om)}~\rho^{n }. $$
\end{enumerate}
\end{theo}

\noindent Revenons maintenant au système initial et essayons d'en déduire la loi invariante associée à $X_n$. Si $(X_n)_n$ est défini par $X_0, X_1$ à valeurs dans $[-L,L]$  et la relation de récurrence $X_{n+2}= \ph(X_n,X_{n+1})$, on pose $Z_n=(X_n,\gamma X_{n+1})$. Alors $(Z_n)_n$ vérifie la relation de récurrence $Z_{n+1}=T(Z_n)$ et on en déduit :

\begin{theo}
Supposons que la variable aléatoire $Z_0=(X_0,\gamma X_{1})$ a pour densité $h_*$. Alors $Z_n$ a pour densité $h_*$. En prenant les marginales, on en déduit que pour tout $n \in \N$, $X_n$ a pour densité
$$
f : x\mapsto \int_{[-\gamma L,\gamma L]} h_*(x,v)\ dv
= \gamma \int_{[- L, L]} h_*(u,\gamma x)\ du.
$$
\end{theo}

\noindent En effet, comme $Z_n= (X_n,\gamma X_{n+1})$ a la densité $h_*$, on obtient que $X_n$ a la densité $f$ en calculant la première marginale. En calculant la deuxième, on a que $\gamma X_{n+1}$ a la densité $g=g(y)$
 définie par
$$
g(y)= \int_{[-L,L]} h_*(u,y)\ du.
$$
On en déduit que $ X_{n+1}$ a la densité $y\mapsto \gamma g(\gamma y)$. Or $Z_{n+1}$ a aussi la densité $h_*$, donc $X_{n+1}$ a aussi la densité $f$ donnée par la première marginale, ce qui donne l'égalité annoncée.\\

\noindent Si $F$ est définie sur $[-L,L]$, notons $Tr ~ F$ la fonction définie, sur $\Omega$, par $Tr ~ F(x,y)= F(x)$.\\
\noindent On a le résultat suivant, conséquence directe du point 6 du théorème 2 appliqué à $Tr ~ F$ et $Tr ~ H$ :

\begin{theo}\label{th4}
Pour tout $B$ borélien et $I$ intervalle, si $(X_0,X_1)$ suit la loi invariante, on a :
$$
\begin{array}{c}
\displaystyle \left\vert P\left( X_{ n \times ppcm(L_i)}\in B, X_0 \in I \right)  -\sum_{j=1}^{\dim(E_1)} \sum_{l=0}^{L_j-1} \mu(W_{j,l}) <Tr ~ \1_B,1>_{\mu_{j,l}} <1,Tr ~ \1_I>_{\mu_{j,l}} \right| \\
\\
\leq 16 (1+\gamma) ~ C ~ L^3 ~ (10 \eps_0^{1-\alpha} + L) ~ \rho^{n }.
\end{array}
$$
Plus généralement, soit $F$ définie, mesurable sur $[-L,L]$ telle que $Tr ~ F$ appartienne à $L^1_{\mu}(\Om)$. Soit $H \in L^{\infty}_m([-L,L])$ telle que $\displaystyle \sup_{0<\eps<\eps_0}\eps^{-\alpha} \int_{[-L,L]}
 {\rm Osc}(H,]x-\eps,x+\eps[ \cap [-L,L])\ dx < +\infty$.\\
 Alors $Tr ~ H \in V_{\alpha}(\Om)$ et 

$$
\left\vert E( F( X_{n \times ppcm(L_i)}) H(X_{0}))
 -\sum_{j=1}^{\dim(E_1)} \sum_{l=0}^{L_j-1} \mu(W_{j,l})
 \mu_{j,l}(Tr ~ F)\mu_{j,l}(Tr ~ H) \right| 
\leq C(F,H) ~ \rho^n $$
avec
$$ \begin{array}{l}
\displaystyle
C(F,H) = C ~ L ~ ||Tr ~ F||_{L^1_{\mu}}\left(2\gamma \sup_{0<\eps<\eps_0}\eps^{-\alpha} \int_{[-L,L]} {\rm Osc}(H,]x-\eps,x+\eps[ \cap [-L,L])\ dx \right. 
\\
\displaystyle + 16 (1+\gamma) ~ \eps_0^{1-\alpha} ||H||_{L^{\infty}_m([-L,L])} + 2\gamma ~ ||H||_{L^1_m([-L,L])}\bigg).
\end{array}
$$
Si de plus $T$ est mélange, alors :
$$ |Cov(F(X_n),H(X_0))| \leq C(F,H) ~ \rho^n. $$
\end{theo}

%
%

\section{Démonstrations}

\noindent Le théorème 2 est une conséquence des théorèmes 5.1 et 6.1 de \cite{SAU}. La difficulté réside dans la vérification, par $T$, des hypothèses (PE1) à (PE5).\\

\noindent  Pour vérifier (PE2), nous allons d'abord montrer que $T_k$ est un $C^1$ difféomorphisme de $W_k$ dans $T_k(W_k)$.
L'hypothèse 3 sur $\displaystyle \frac{\partial \ph_k}{\partial u}$ assure que $T_k$ est un difféomorphisme local.
Pour l'injectivité, considérons $(x,y)$ et $(x',y')$ deux points différents de $W_k$, ayant même image par $T_k$. On obtient $y=y'$ et $\ph_k(x',y/\gamma)= \ph_k(x,y/\gamma)$.
En utilisant l'hypothèse géométrique 4 et en appliquant le théorème des accroissements finis à $t \mapsto \ph_k(\Gamma_1(t),\Gamma_2(t) )$, on aboutit à une contradiction.

\noindent Les hypothèses de régularité sur les $\ph_k$ (et donc les $T_k$) permettent d'établir que $\det(DT_k^{-1})$ est höldérien d'exposant $\alpha$, si l'on restreint convenablement le domaine : on peut voir qu'il existe, pour chaque $k$, $\beta_{k}>0$, un ouvert $\Val_k$ d'adhérence compacte et une constante $c_{k}$ tels que 
\begin{itemize}
\item $\overline{U_k}\subset  \Val_{k} \subset \overline{\Val_{k}} \subset W_k$ ;
\item $B_{\beta_{k}}(T_k(U_k)) \subset T_k(\Val_{k}) ;$
\item
pour tout $\eps<\beta_{k}$, tout $z\in T_k(\Val_{k})$ et tous $x,y\in B_{\eps}(z)\cap T_k(\Val_{k})$, on ait
$$
\Big| \det(DT_k^{-1}(x))- \det(DT_k^{-1}(y)) \Big| \leq c_{k} \Big| \det(DT_k^{-1}(z)) \Big|\eps^{\alpha}.
$$
\end{itemize}

\noindent En posant $\beta = \min\limits_k \beta_{k} >0$ et $c = \max\limits_k c_{k} >0$, on obtient des constantes valables pour tout $k \in \{1, \hdots , d\}$, d'où (PE2) est vérifiée.
\\

\noindent Cela permet de fixer les ouverts sur lesquels on va travailler : il existe $\eps_2>0$ tel que, pour tout $k\in\{1,\dots,d\}$, $B_{2\eps_2}(\overline{U_k}) \subset \Val_{k} \subset W_k$.
On prend désormais $V_k=B_{\eps_2}(\overline{U_k})$. On a $T_k(V_k)$ ouvert et $T_k(\overline{U_k})$ compact inclus dans $T_k(V_k)$. On peut trouver un $\eps_0^1 >0$ tel que $B_{\eps_0^1}( T_k(\overline{U_k})) \subset T_k(V_k)$ pour tout $k$. L'hypothèse (PE1) est ainsi vérifiée.\\
\\
L'hypothèse (PE3) est évidente car $\displaystyle \Om= \bigcup_{k=1}^d  U_k\ \cup \Nal'$ est une réunion disjointe d'ouverts et d'une partie négligeable.\\
\\
Pour (PE4), on procède en deux étapes : d'abord montrer un résultat de dilatance lorsque les points antécédents dans $\Val_k$ sont proches (proposition \ref{dilatance1}) puis l'hypothèse proprement dite (proposition \ref{dilatance2}), résultat de dilatance lorsque les points images dans $T_k(V_k)$ sont proches.

\begin{prop}\label{dilatance1}
Soient $(x,y)$ et $(x',y')\in \Val_{k}$ tels que le segment $[(x,y),(x',y')]$ soit inclus dans $\Val_{k}$. Alors 
$$
|| T_k(x,y)-T_k(x',y')||^2 \geq \frac{1}{s^2} || (x,y)-(x',y')||^2.
$$
\end{prop}

\noindent{\it Démonstration :} On applique le théorème des accroissement finis à l'application définie sur $[0,1]$ par $t \mapsto \ph_k(x+t(x'-x), \frac{1}{\gamma} (y+t(y'-y))$, ce qui nous permet de voir qu'il existe $c\in ]0,1[$ tel que 
$$
 || T_k(x,y)-T_k(x',y')||^2 = (x'-x,y'-y) B \left( \begin{array}{lll}
x'-x\\
y'-y
\end{array}\right)
$$
où
$$
B = \left( \begin{array}{lll}
\displaystyle
\gamma^2 
\left(\frac{\partial \ph_k}{\partial u}(x_c,\frac{1}{\gamma}y_c)\right)^2
&
\displaystyle
\gamma \frac{\partial \ph_k}{\partial u}(x_c,\frac{1}{\gamma}y_c)\frac{\partial \ph_k}{\partial v}(x_c,\frac{1}{\gamma}y_c)\\
\displaystyle
\gamma \frac{\partial \ph_k}{\partial u}(x_c,\frac{1}{\gamma}y_c)\frac{\partial \ph_k}{\partial v}(x_c,\frac{1}{\gamma}y_c)  &
\displaystyle \frac{1}{\gamma^2} +
\left(\frac{\partial \ph_k}{\partial v}(x_c,\frac{1}{\gamma}y_c)\right)^2
\end{array}\right)
$$
avec $(x_c,y_c) = (x+c(x'-x),y+c(y'-y))$.\\
La matrice $B$ est symétrique, réelle. Posons 
$$
\begin{array}{lll}
\displaystyle
\xi_1 & = {\rm Tr}(B) & \displaystyle = 
 \frac{1}{\gamma^2} +
\left(\frac{\partial \ph_k}{\partial v}(x_c,\frac{1}{\gamma}y_c)\right)^2
+ \gamma^2 
\left(\frac{\partial \ph_k}{\partial u}(x_c,\frac{1}{\gamma}y_c)\right)^2 \\
\xi_2 & =\det(B) &\displaystyle = \left(\frac{\partial \ph_k}{\partial u}
(x_c,\frac{1}{\gamma}y_c)\right)^2.
\end{array}
$$
Nous allons établir que les valeurs propres de $B$ sont supérieures ou égales à $\frac{1}{s^2}$. En effet, l'application $\zeta :\R^2 \rightarrow \R^2$ définie par $\zeta(x,y)= (x+y,xy)$ réalise une bijection entre les ensembles 
$$
V''_s=\{ (x,y) \in \R^2 : s^{-2} \leq x\leq y\}\quad {\rm et }\quad \zeta(V''_s)= \{ (\xi_1,\xi_2) \in (\R^*_+)^2 : \xi_1 \geq 2s^{-2},\ \xi_2\geq s^{-2}(\xi_1-s^{-2}),\ \xi_2 \leq \frac{\xi_1^2}{4} \}.
$$
Il suffit donc de vérifier que $(\xi_1,\xi_2)$ est dans $\zeta(V''_s)$ pour avoir le résultat.\\
Or, les valeurs propres de $B$ sont réelles, le discriminant de son polynôme caractéristique est positif donc $4\xi_2\leq \xi_1^2$ et les conditions sur $A$ et $M$, ainsi que le choix de $s$ et $\gamma$, assurent que les deux autres inégalités sont vérifiées.\\
Donc la matrice $B$ a ses valeurs propres supérieures ou égales à $s^{-2}$. Il s'ensuit que $ || T_k(x,y)-T_k(x',y')||^2 \geq \frac{1}{s^2} || (x,y)-x',y')||^2$.\\

\noindent Par des arguments de compacité, on montre qu'il existe $\eps_0^2 > 0$ tel que, pour tout $(x,y)\in \overline{V_k}$,
$$
B_{\eps_0^2}(T_k(x,y)) \subset T_k(B_{\eps_2}(x,y)).
$$

\begin{prop}\label{dilatance2}
Soit $\eps_0= \min(\eps_0^1,\eps_0^2) > 0$. On rappelle que $\overline{U_k} \subset V_k \subset \overline{V_k} \subset \Val_k \subset W_k$. Alors :
\begin{itemize}
\item
Quels que soient $(u_1,v_1), (u_2,v_2) \in T_k(V_k)$ vérifiant $d((u_1,v_1), (u_2,v_2))<\eps_0$, on a l'inégalité 
$$
s^2 ~ d((u_1,v_1), (u_2,v_2)) >
 d(T_k^{-1}(u_1,v_1),(T_k^{-1} (u_2,v_2)).
$$
avec $\displaystyle s^2 < 1$.
\item
$B_{\eps_0}( T_k(\overline{U_k})) \subset T_k(V_k)$.
\end{itemize}
\end{prop}

\noindent{\it Démonstration :} La deuxième affirmation vient de ce que $\eps_0\leq \eps_0^1$ et de ce qu'on a obtenu dans (PE1). \\
La première affirmation entraîne la condition (PE4) de Saussol. En effet, soient $(u_1,v_1), (u_2,v_2) \in T_k(V_k)$ vérifiant $d((u_1,v_1), (u_2,v_2))<\eps_0$. Soit $(x,y)=T_k^{-1}(u_1,v_1) \in V_k$. D'après la remarque précédente, $\eps_0$ étant inférieur à $\eps_0^2$,
$$
(u_2,v_2) \in B_{\eps_0}(T_k(x,y)) \subset T_k(B_{\eps_2}(x,y)).
$$
Donc $(x',y')=T_k^{-1}(u_2,v_2)\in B_{\eps_2}(x,y)\subset \Val_k$. D'après la proposition \ref{dilatance1}, 
$$
d((u_1,v_1),(u_2,v_2))^2 = 
|| T_k(x,y)-T_k(x',y')||^2 > \sigma || (x,y)-(x',y')||^2,
$$
ce qui prouve le résultat.\\

\noindent Pour finir, l'hypothèse (PE5) vient du lemme 2.1 de Saussol et de l'hypothèse 6.\\

\noindent Puisque les hypothèses (PE1) à (PE5) sont vérifiées, le théorème 5.1 de \cite{SAU} implique les propriétés 1 à 5 du théorème 2 sur $V_{\alpha}$ et $L^1_m$. Mais, si $f \in E_i$, $f$ est nulle sur $\Omega^c$ et donc $f \in L^1_m(\Omega)$ et $V_{\alpha}(\Omega)$.\\

\noindent Pour démontrer le point 6, nous allons appliquer le théorème 6.1 de \cite{SAU}
sur chaque $W_{j,l}$ où une puissance de $T$ est mélange. Avec les notations du
 point 5 du théorème 5.2 de  \cite{SAU}, on a l'existence de constantes $C>0$
et $\rho\in ]0,1[$ telles que, pour tout $(j,l)$ vérifiant
$ 1\leq j\leq \dim(E_1)$, $0\leq l\leq L_j -1$, toute fonction 
$f\in L^1_{\mu_{j,l}}(\Om)$ et pour toute fonction $h\in V_{\alpha}(\Om)$,
$$
\left\vert \int_{\Om}(f- \mu_{j,l}(f))\circ T^{nL_j} h ~ d\mu_{j,l}
\right\vert \leq C ||f-\mu_{j,l}(f))||_{L^1_{\mu_{j,l}}} ||h||_{\alpha,L}\rho^{n } .
$$
Soient alors $h\in V_{\alpha}(\Om)$ et  $f\in L^1_{\mu}(\Om)$ (de sorte que
$f\in L^1_{\mu_{j,l}}(\Om)$ pour chaque $j,l$). En prenant le p.p.c.m. $L'$ des $L_j$ et en
faisant la somme des inégalités ci-dessus, avec $n$ remplacé par
$n \frac{L'}{L_j}$, on obtient
$$
\left| \int_{\Om} f \circ T^{ n L'}  h \ d\mu -\sum_{j=1}^{\dim(E_1)}
 \sum_{l=0}^{L_j-1} \mu(W_{j,l}) \mu_{j,l}(f)\mu_{j,l}(h) \right| 
\leq C  ||h||_{\alpha,\Om}||f||_{L^1_{\mu}} \rho^{n } .
$$

\noindent Le point 7 est une conséquence directe du 6, vu que dim($E_1) = 1$ et que $L_1=1$.\\

\noindent Passons au théorème 4. Si $\left( \begin{array}{lll} X_{0}\\ \gamma X_{1} \end{array}\right)$ suit la loi $\mu$, alors  $\left( \begin{array}{lll} X_{n}\\ \gamma X_{n+1} \end{array}\right)$ aussi. Si $f\in L^1_{\mu}(\Om)$ et si $h\in V_{\alpha}(\Omega)$, on a :
$$
\left\vert E\left( f \left( \begin{array}{lll} X_{nL'}\\ \gamma X_{nL'+1} \end{array}\right) h\left( \begin{array}{lll} X_{0}\\ \gamma X_{1} \end{array}\right)\right)
 -\sum_{j=1}^{\dim(E_1)} \sum_{l=0}^{L_j-1} \mu(W_{j,l}) \mu_{j,l}(f)\mu_{j,l}(h) \right| \leq C ||f||_{L^1_{\mu}} ||h||_{\alpha,\Om}\rho^{n } .
$$
\bigskip

\noindent Pour que $Tr~ H$ appartienne à $V_{\alpha}(\Omega)$, il faut et il suffit que 
$H$ soit dans $L^{\infty}([-L,L],m)$ et vérifie
$$
 \sup_{0<\eps<\eps_0}\eps^{-\alpha} 
\int_{[-L,L]} {\rm Osc}(H,]x-\eps,x+\eps[\cap [-L,L])
\ dx  <\infty.
$$

 De plus,
$$
\begin{array}{lll}
|| Tr ~ H||_{\alpha,\Om} 
 & \displaystyle =
2\gamma L \sup_{0<\eps<\eps_0}\eps^{-\alpha} \int_{[-L,L]}
 {\rm Osc}(H,]x-\eps,x+\eps[ \cap [-L,L])\ dx\\
& \displaystyle
 + 16(1+\gamma) L \eps_0^{1-\alpha} ||H||_{L^{\infty}_m([-L,L])}
 + 2\gamma L ||H||_{L^1_m([-L,L])}.
\end{array}
$$
Donc si $H$ vérifie ces conditions et si $F $ est telle que 
$Tr~F$ soit dans $L^1_{\mu}(\Om)$, par exemple si $F$ est mesurable et bornée sur
$[-L,L]$, on a

$$
\begin{array}{lllll}
\displaystyle
\left\vert E( F( X_{n \times L'}) H(X_{0}))
 -\sum_{j=1}^{\dim(E_1)} \sum_{l=0}^{L_j-1} \mu(W_{j,l})
 \mu_{j,l}(Tr ~ F)\mu_{j,l}(Tr ~ H) \right| 
\\
\displaystyle\leq C 
 ||Tr ~ F||_{L^1_{\mu}}\left(
2\gamma L \sup_{0<\eps<\eps_0}\eps^{-\alpha} \int_{[-L,L]}
 {\rm Osc}(H,]x-\eps,x+\eps[ \cap [-L,L])\ dx \right. 
\\
\displaystyle + 16 (1+\gamma) L \eps_0^{1-\alpha} ||H||_{L^{\infty}_m([-L,L])} + 2\gamma L ||H||_{L^1_m([-L,L])}\bigg)\rho^{n } .
\end{array}
$$

\noindent En particulier, si $H$ est l'indicatrice d'un intervalle et $F$, celle d'un borélien, on obtient la première assertion du théorème \ref{th4}.

%
%

\section{Exemples}

\subsection{Un cas non linéaire}

\noindent On note, pour tout $k \in \Z$, $f_k$ la fonction polynomiale $\displaystyle f_k(x) = - \frac{71}{2} x^2 - 214 x + k - \frac{1}{2}$.\\
Pour tout $ -179 \leq k \leq 250$, on définit l'ouvert $O_k$ par :
$$ O_k = \{ (u,v) \in ]-1,1[^2 ~ / ~ f_k(u) < v < f_{k+1} (u) \}. $$

\noindent On considère les applications $\ph_k$ définies sur $B_1(\overline{O_k})$ ($\eps_1 = 1$) pour tout $-179 \leq k \leq 250$ par :
$$ \ph_k(u,v) = 2v - 2 f_k(u) - 1. $$
On définit $\ph : [-1,1]^2 \rightarrow [-1,1]$ presque partout en posant $\ph_{|_{O_k}} = \ph_{k|_{O_k}}$ pour tout $-179 \leq k \leq 250$.
Nous allons maintenant nous assurer que ces fonctions et ouverts vérifient les conditions de la partie 2.\\

\noindent La condition sur les ouverts est facilement vérifiée, puisque $[-1,1]^2 \backslash \bigcup\limits_{k=-179}^{250} O_k $ est une réunion de segments et de paraboles. De plus, le nombre maximal d'arcs se croisant est $Y=3$. \\
Les conditions de régularité sont satisfaites car les $\ph_k$ sont $C^{\infty}$ sur $B_1(\overline{O_k})$. On prend $\alpha = 1$. Les dérivées partielles vérifient les inégalités suivantes : pour tout $-179 \leq k \leq 250$ et pour tout $(u,v) \in B_1(\overline{O_k})$ on a :
$$ \left| \frac{\partial \ph_k}{\partial v} (u,v) \right| = 2 = M $$
et 
$$ \left| \frac{\partial \ph_k}{\partial u} (u,v) \right| = 2 | 71 u + 214| > 2 (214 - 71 (1 + 1)) = 144 = A > M+1. $$

\noindent Dans ce cas, $\gamma = \frac{1}{12}$. Un calcul montre que $s \leq \frac{1}{10}$ et $\eta < 1$.\\
On pose $\Om= [-1,1]\times\left[ - \frac{1}{12},  \frac{1}{12}\right]$ et pour tout $-179 \leq k \leq 250$ on définit les ouverts 
$$
U_k=  \{(x,y)\in \mathring{\Om} : f_k(x) < 12 y < f_{k+1}(x) \}.
$$
On obtient les applications 
$$
T_k(x,y)= (12 y, \frac{2}{12}(12 y -f_k(x))- \frac{1}{12}).
$$

\noindent Si $-177 \leq k\leq 248$, $T_k(U_k)= \mathring{\Om}$ et $T_k$ est bijectif de $U_k$ sur $\mathring{\Om}$.\\
Sinon, on peut vérifier que $T_{-178}$ est une bijection de $U_{-178}$ sur $\displaystyle \Om_1 \cup \Om_2$, où $\Omega_1$ est la partie ouverte de $\mathring{\Omega}$ située au-dessus de la droite d'équation $y= \frac{2x+1}{12}$, $\Omega_3$ celle située sous la droite $y=\frac{2x-1}{12}$ et $\Omega_2$ celle située entre ces deux droites. On a des relations similaires pour $k=-177$, $249$ et $250$ et d'autres sous-ensembles de $\Om$.\\

\noindent Enfin la condition géométrique est vérifiée sous sa forme simple (l'ouvert contient le segment horizontal).\\
La transformation $T$ admet donc une densité invariante $h_*$.\\

\noindent Soit $P$ l'opérateur de Perron Frobenius associé à $T$. On peut vérifier que les fonctions constantes ne sont pas invariantes par $P$, donc que $h_*$ n'est pas constante.
On pose 
$$
\psi_k(x,y)= (214)^2-71(2x-12y)+142 k.
$$
Alors $Ph$ a l'expression $\displaystyle Ph(x,y)= \sum_{k=a}^{b} h(T_k^{-1}(x,y)) \frac{1}{2 \sqrt{\Psi_k(x,y)}}$, avec $(a,b)=(-179,248)$ si $(x,y) \in \Omega_1$, $(a,b) = (-178,249)$ si $(x,y) \in \Omega_2$, $(a,b) = (-177,250)$ si $(x,y) \in \Omega_3$.\\

\noindent On va établir que $P1\neq 1$. Supposons que $h=1$ et posons $z=x-6y$. \\
Si $(x,y)\in \Om_3$, $z \in ]-\frac{3}{2}, -\frac{1}{2}[$. La fonction $\displaystyle z\mapsto \sqrt{ (214)^2-142 z +142 k}$
est strictement décroissante sur $]-\frac{3}{2}, -\frac{1}{2}[$. Donc $\displaystyle z\mapsto  \frac{1}{2 \sqrt{ (214)^2-71(2x-12y)+142 k}}$ est strictement croissante sur $]-\frac{3}{2}, -\frac{1}{2}[$ et
$P1$ n'est donc pas constante.

\subsection{Un cas linéaire par morceaux}

\noindent Dans toute cette partie, $a$ et $b$ sont des entiers relatifs non nuls, $L$ est un entier ou demi entier strictement positif.\\
On note $\Ual^2$ le carré $[-L,L]^2$. Pour tout $n\in \Z$,  l'ouvert $\Omega_n$ est défini par
$$
\Omega_n=\{ (u,v)\in ]-L,L[^2 \ : \  av+bu \in ] (2n-1)L,(2n+1)L[\}.
$$
On désigne par $\Delta_n$ la droite d'équation $ av+bu=(2n-1)L$. On définit $\varphi_n$ sur $\R^2$ par 
$$
\varphi_n(u,v)=  av+bu-2nL.
$$
Alors  $\left.\varphi_n\right|_{\Omega_n}$ est à valeurs dans $]-L,L[$ et l'on pose
$$
\forall (u,v)\in \Omega_n,\ \varphi(u,v) = \varphi_n(u,v).
$$
On impose la condition suivante, avec $\displaystyle S= 1+ \frac{48}{\pi} + \frac{288}{\pi^2} + 
\frac{4}{\pi} \left(1+\frac{12}{\pi}  \right)\sqrt{6\pi +36}$,

$$ |a|< \frac{ |b| -S}{\sqrt{S}}. $$
\vskip 0.5 cm

\noindent On vérifie les conditions de la partie 2.\\
Le carré $\Ual^2$ est la réunion disjointe d'un nombre fini d'ouverts $\Om_n$ et d'une partie négligeable constituée d'une nombre fini de segments de droites.\\
Le nombre maximal de croisements de ces segments est $Y=3$.\\
Pour tout $\eta>0$, les ouverts $ B_{\eta}(\Omega_n)$ sont convexes donc le chemin $\Gamma$ joignant deux points de même ordonnée est le segment horizontal et la condition géométrique est satisfaite.\\
Les $\varphi_n$ sont de classe $C^{\infty}$ sur $ B_{\eta}(\Omega_n)$. On prendra $\alpha=1$. De plus  $\varphi_n(\Omega_n)\subset [-L,L]$.\\
On pose $M=|a|$, $A=|b|$, de sorte que les inégalités sur les dérivées partielles soient satisfaites. La majoration de $|a|$ implique que $0<M<A-1$.\\
On définit $\gamma= |b|^{-1/2}$, le coefficient d'aplatissement et on vérifie que $\eta<1$.\\

\noindent On détermine pour quels entiers $n$ la droite $\Delta_n$ intersecte le carré $\Ual^2$. On vérifie que $\Delta_n\cap \Ual^2 \neq \emptyset$ si et seulement si
$$
\frac{-|a|- |b|+1}{2}\leq n\leq \frac{|a|+|b|+1}{2}.
$$

\noindent On définit $\nnn(a,b)$ comme étant l'ensemble des indices $n$ tels que $\Omega_n \neq \emptyset$  intersecte $\Ual^2$.

\noindent Posons
$$
\Omega= [-L,L]\times [-\gamma L,\gamma L]
$$
et 
$$
\Omega_{n,a}=
\{ (x,y)\in \mathring{\Omega}\ :\ a \sqrt{|b|} y + b x -2nL \in ]-L,L[ \}
$$
On pose 
$$
\begin{array}{lllll}
T_n & : & \Omega_{n,a} & \rightarrow & \Omega\\
&& (x,y) & \mapsto &  ( \sqrt{|b|} y, ay+ \frac{b}{ \sqrt{|b|} }x
- \frac{2nL}{ \sqrt{|b|} } ). \\
\end{array}
$$
On définit $T$ presque partout, de $\Omega$ dans lui-même en posant $\left. T \right|_{\Omega_{n,a}} = T_n$.\\

\noindent L'opérateur de Perron Frobenius $P$ associé à $T$ a l'expression 
$$
Ph(x,y)=\frac{1}{|b|}
\sum_{n\in \nnn(a,b)} \1_{(x,y)\in T_n(\Omega_{n,a}) } h(T_n^{-1}(x,y)).
$$
Si $h$ est une fonction constante égale à $c>0$, on en déduit que 
$$
Ph(x,y)= \frac{c}{|b|} \sharp \{ n\in \nnn(a,b) \ : \ (x,y)\in 
 T_n(\Omega_{n,a})  \}.
$$
On vérifie que ce cardinal est bien $|b|$, ce qui donne une densité invariante constante.

\noindent Or le théorème donne l'existence d'une mesure invariante $h^*m$, donnée par $P_1 \1_{\Om} = h^*$. D'après le lemme 4.1 de \cite{ITM}, 
$$ P_1 \1_{\Om} =\lim\limits_{n \rightarrow +\infty} \frac{1}{n} \sum\limits_{k=1}^n P^k \1_{\Om} = \1_{\Om}.
$$

%
%


\begin{thebibliography}{99}

\bibitem[AFLV]{AFLV}
 ALVES José F.,  FREITAS Jorge M., LUZZATTO Stefano, VAIENTI Sandro,
{\it  From rates of mixing to recurrence times via large deviations,  Advances in Mathematics},  \textbf{228} (2011), n° \textbf{2} 1203-1236.
\bibitem[CE]{CE}
COLLET Pierre, ECKMANN  Jean-Pierre, {\it Concepts and results in chaotic dynamics: a short course}. Theoretical and Mathematical Physics. Springer-Verlag, Berlin (2006).
\bibitem[HK]{HK}
HOFBAUER Franz, KELLER Gerhard, {\it Ergodic properties of invariant measures for piecewise monotonic transformations}, Mathematische Zeitschrift {\bf 180} (1982), 119-140.
\bibitem[ITM]{ITM}
IONESCU TULCEA C.T., MARINESCU G., {\it Théorie ergodique pour des classes d'opérations non complètement continues}, Annals of Mathematics {\bf Vol. 52, n°2} (1950), 140-147.
\bibitem[LM]{LM}
LASOTA Andrzej, MACKEY Michael C., {\it Chaos, fractals and noise : stochastic aspects of dynamics}, Springer Verlag, New York (1998)
\bibitem[LIV]{LIV}
LIVERANI Carlangelo, {\it Multidimensional expanding maps with singularities: a pedestrian approach}, Ergodic Theory and Dynamical Systems {\bf Vol. 33, n°1} (2013), 168-182.
\bibitem[SAU]{SAU}
SAUSSOL Benoît, {\it Absolutely continuous invariant measures for multidimensional expanding maps}, Israel Journal of Mathematics {\bf 116} (2000), 223-248.
\bibitem[TON1]{TON1}
 TONG Howell, {\it Nonlinear time series. A dynamical system approach}.
  With an appendix by K. S. Chan. Oxford Statistical Science Series, \textbf{6}. Oxford Science Publications. The Clarendon Press, Oxford University Press, New York  (1990).

\bibitem[TON2]{TON2}
    TONG Howel,  {\it Nonlinear time series analysis since 1990: some personal reflections}. Acta Math. Appl. Sin. Engl. Ser. 18 (2002), no. \textbf{2}, 177-184.

\bibitem[YOU]{YOU}
YOUNG Lai-Sang, {\it Recurrence times and rates of mixing}, Israel Journal of Mathematics {\bf 110} (1999), 153-188.

\end{thebibliography}
\end{document}